\input xy
\xyoption{all}\xyoption{poly}
\newdir{ >}{{}*!/-9pt/@{>}}

\documentclass[reqno]{amsart}
\usepackage{amsmath,latexsym,amssymb,mathrsfs}
\usepackage{hyperref}
\usepackage{textcomp}
\usepackage{graphicx}

\usepackage{hyperref}
\usepackage[all]{xy}

\theoremstyle{plain}
\newtheorem{theorem}{Theorem}[section]
\newtheorem{lemma}[theorem]{Lemma}
\newtheorem{proposition}[theorem]{Proposition}

\theoremstyle{definition}

\newtheorem{examples}[theorem]{Examples}

\theoremstyle{remark}
\newtheorem{remark}[theorem]{Remark}

\newcommand{\CC}{ \ensuremath{\mathbb {C}} }

\newcommand{\map}[2]{ \ensuremath{ \xymatrix@1@C=15pt{ #1 \ar[r] & #2 } } }
\newcommand{\mono}[2]{ \ensuremath{ \xymatrix@1@C=15pt{ #1 \ar@{ >->}[r] & #2 } } }
\newcommand{\regepi}[2]{ \ensuremath{ \xymatrix@1@C=15pt{ #1 \ar@{>>}[r] & #2 } } }
\renewcommand{\AA}{ \ensuremath{\mathbb {A}} }
\newcommand{\EE}{ \ensuremath{\mathcal {E}} }
\newcommand{\Grp}{ \ensuremath{\mathrm{Grp}} }

\begin{document}

\title{The Cuboid Lemma and Mal'tsev categories}
\thanks{Research supported by the F.N.R.S. grant {\em Cr\'edit aux chercheurs} 1.5.016.10F, and by the Centro de Matem\'atica da
Universidade de Coimbra (CMUC), funded by the European Regional
Development Fund through the program COMPETE and by the Portuguese
Government through the FCT - Funda\c{c}\~ao para a Ci\^encia e a Tecnologia
under the project PEst-C/MAT/UI0324/2011.}

\author{Marino Gran     \and
        Diana Rodelo 
}


\address{Marino Gran \\
             Insitut de Recherche en Math\'ematique et Physique, Universit\'e catholique de Louvain, Chemin du Cyclotron 2,
1348 Louvain-la-Neuve, Belgium}
    \address{
           Diana Rodelo \\
           Departamento de Matem\'atica, Faculdade de Ci\^{e}ncias e Tecnologia, Universidade do Algarve, Campus de
Gambelas, 8005--139 Faro, Portugal\newline \and  CMUC, Universidade de Coimbra, 3001--454 Coimbra, Portugal
}
 \email{marino.gran@uclouvain.be \\ drodelo@ualg.pt} 
\maketitle
{\rm{\emph{ Dedicated to George Janelidze on the occasion of his 60th birthday}}}

\begin{abstract}
We prove that a regular category $\CC$ is a Mal'tsev category if and only if a strong form of the denormalised $3 \times 3$ Lemma holds true in $\CC$. In this version of the $3 \times 3$ Lemma, the vertical exact forks are replaced by pullbacks of regular epimorphisms along arbitrary morphisms. The shape of the diagram it determines suggests to call it the Cuboid Lemma. This new characterisation of regular categories that are Mal'tsev categories (= 2-permutable) is similar to the one previously obtained for Goursat categories (= 3-permutable). We also analyse the ``relative'' version of the Cuboid Lemma and extend our results to that context. \\

\noindent Math. Subj. Class. 2010: 18C05, 08C05, 18B10, 18E10. 
\end{abstract}

\section{Introduction}
\label{intro}
One of the motivations to introduce the notion of a \emph{Mal'tsev category}, already mentioned in the
 classical article by A. Carboni, J. Lambek and M.C. Pedicchio \cite{CLP}, was the possibility of establishing extended versions of homological diagram lemmas, such as
the Snake Lemma, in a non-pointed context.
This idea was further pursued in the work of D. Bourn \cite{B}, who established an extended version of the  $3 \times 3$ Lemma in any regular Mal'tsev category, called the \emph{denormalised  $3 \times 3$ Lemma}. In general, in order to formulate this kind of result in a non-pointed context, the short exact sequences
are simply replaced by \emph{exact forks}, i.e. by diagrams of the form
$$
\xymatrix{
  R_f \ar@<3pt>[r]^-{f_1} \ar@<-3pt>[r]_-{f_2} & A \ar@{>>}[r]^-f & B }
$$
where $f$ is a regular epimorphism, and $(R_f,f_1,f_2)$ is the kernel pair of $f$.

The validity of the denormalised  $3 \times 3$ Lemma is actually equivalent to the property of $3$-permutability of the composition of equivalence relations \cite{L,GR}: for all equivalence relations $R$ and $S$ on the same object one has that $$R S R =SRS.$$
Those regular categories having $3$-permutable equivalence relations are called \emph{Goursat categories}, and they were introduced in \cite{CKP}.
However, the regular categories having the Mal'tsev property are precisely those that satisfy the stronger property of $2$-permutability of the composition of equivalence relations, namely
 $$R S  =SR.$$

These observations naturally led to the problem of determining whether there existed a homological diagram lemma whose validity would characterise Mal'tsev categories among regular ones, in the same way as the denormalised  $3 \times 3$ Lemma characterises Goursat categories. The present article answers this question: the \emph{Cuboid Lemma}
 is a new diagram lemma, stronger than the denormalised  $3 \times 3$ Lemma, allowing one to obtain the desired characterisation of Mal'tsev categories (Theorems \ref{thm upper split cuboid} and \ref{thm upper cuboid}). This result is obtained by using the ``calculus of relations'' available in any regular category, much in the spirit of the above mentioned article \cite{CLP}.

 Our new observations can also be extended from the ``absolute'' context of regular categories to the wider context of ``relative'' regular categories~\cite{TJ,RS-AC,RGC}. Here, the role of regular epimorphisms (for regular categories) is played by morphisms that belong to a distinguished class $\EE$ of regular epimorphisms. Relative regular categories give a suitable setting to develop a ``calculus of $\EE$-relations'': this makes it possible to extend the main results obtained in a regular category to the ``relative'' context, by using exactly the same proofs.

\vspace{3mm}

\noindent {\bf Structure of the paper}\\
In Section~\ref{Preliminaries} we give the main background concerning regular categories and the calculus of relations. We revise some known properties of regular Mal'tsev categories and their characterisation through the so-called regular pushouts in Section~\ref{Regular Mal'tsev categories}. We also include a new characterisation of Mal'tsev categories through a suitable stability property of regular epimorphisms. In Section~\ref{The Cuboid Lemma} we state the Cuboid Lemma and prove that it characterises regular categories that are Mal'tsev categories. We conclude the article by extending the main results of this paper to a relative context in Section~\ref{The relative context}.

\section{Preliminaries}\label{Preliminaries}
A finitely complete category $\CC$ is said to be a \emph{regular} category when any kernel pair has a coequaliser and, moreover,
regular epimorphisms are stable under pullbacks. In a regular category any morphism $f \colon A \rightarrow B$ has a factorisation $f=m\cdot p$,  where $p$  is a regular epimorphism and $m$ is a monomorphism. The corresponding (regular epimorphism, monomorphism) factorisation system is then stable under pullbacks.

\noindent \emph{In this article $\CC$ will always denote a finitely complete regular category.}

 A relation $R$ from $A$ to $B$ is a subobject $\langle r_1,r_2 \rangle\colon R\rightarrowtail A\times B$. The opposite relation, denoted $R^{\circ}$, is the relation from $B$ to $A$ given by the subobject $\langle r_2,r_1 \rangle\colon R\rightarrowtail B\times A$. We identify a morphism $f:A\rightarrow B$ with the relation $\langle 1_A,f \rangle\colon A\rightarrowtail A\times B$ and write $f^{\circ}$ for the opposite relation. Given another relation $S$ from $B$ to $C$, the composite relation of $R$ and $S$ is a relation $SR$ from $A$ to $C$. With the above notation, we can write any relation $\langle r_1,r_2 \rangle\colon R\rightarrowtail A\times B$ as $R=r_2 r_1^{\circ}$. The following properties are well known (see \cite{CKP}, for instance); we collect them in a lemma for future references.
\\
\begin{lemma}
\label{pps of ms as relations}
Let $f: A\rightarrow B$ be any morphism in a regular category $\CC$. Then:
\begin{enumerate}
 \item {$f f^{\circ} f=f$ and $f^{\circ} f f^{\circ}=f^{\circ}$;}
  \item{
  $f f^{\circ}=1_B$ { if and only if} $f$
  is a regular {epimorphism}.}
\end{enumerate}
\end{lemma}

A relation $R$ from an object $A$ to $A$ is called a \emph{relation on $A$}. Such a relation is \emph{reflexive} if $1_A\leqslant R$, \emph{symmetric} if $R^{\circ} \leqslant R$, and \emph{transitive} when
$RR \leqslant R$. As usual, a relation $R$ on $A$ is an \emph{equivalence relation}  when it is reflexive, symmetric and transitive. In particular, a kernel pair $\langle f_1,f_2 \rangle\colon R_f \rightarrowtail A\times A$ of a morphism $f: A\rightarrow B$ is an \emph{effective equivalence relation}, which can be written either as $R_f=f^{\circ}f$, or as $R_f=f_2f_1^{\circ}$, as mentioned above. When $f$ is a regular epimorphism, then $f$ is the coequaliser of
$f_1$ and $f_2$ and the diagram
$$
\xymatrix{
  R_f \ar@<3pt>[r]^-{f_1} \ar@<-3pt>[r]_-{f_2} & A \ar@{>>}[r]^-f & B }
$$
is called an \emph{exact fork}. Note that, if $f=m\cdot p$ is the (regular epimorphism, monomorphism) factorisation of an arbitrary morphism $f$, then $R_f$=$R_p$, so that an effective equivalence relation is always the kernel pair of a regular epimorphism.

\section{Regular Mal'tsev categories}
\label{Regular Mal'tsev categories}
\noindent Recall that a finitely complete category $\CC$ is called a \emph{Mal'tsev category} when any reflexive relation in $\CC$ is an equivalence relation \cite{CLP,CKP}. The following well known characterisation of regular Mal'tsev categories will be useful.

\begin{proposition}
\label{Definition: Mal'tsev category}
A regular category $\CC$ is a Mal'tsev category if and only if the composition of effective equivalence relations on any object in $\CC$ is commutative:
$$R_f R_g=R_g R_f$$ for any pair of regular epimorphisms $f$ and $g$ in $\CC$ with the same domain.
\end{proposition}

\begin{examples}
{\rm A variety of universal algebras is a Mal'tsev category if and only if its theory has a ternary operation $p(x,y,z)$ satisfying the identities $p(x,y,y)=x$ and $p(x,x,y)=y$ \cite{S}. Of course, the variety $\Grp$ of groups is a Mal'tsev category, as is any variety whose theory contains a binary operation satisfying the group identities. Also the variety of quasi-groups and the variety of Heyting algebras are Mal'tsev categories, as is the dual category of an elementary topos. The category of topological groups is a regular Mal'tsev category, as is the category of Hausdorff groups \cite{CKP}. If $\CC$ is a finitely complete category, then the category $\Grp(\CC)$ of internal groups in $\CC$ is a Mal'tsev category.}
\end{examples}

There are many known characterisations of Mal'tsev categories (see \cite{BB}, for instance, and references therein). We shall now focus on the relationship with commutative diagrams of the form
\begin{equation}
\label{regular po}
\xymatrix{
  C \ar@{>>}[r]^-{c} \ar@<-2pt>[d]_-g & A  \ar@<-2pt>[d]_-f \\
  D \ar@{>>}[r]_-d \ar@<-2pt>[u]_-t & B, \ar@<-2pt>[u]_-s }
\end{equation}
where $f$ and $g$ are split epimorphisms ($f\cdot s=1_B$, $g\cdot t=1_D$), $f\cdot c= d\cdot g$, $s\cdot d= c\cdot t$, and $c$ and $d$ are
regular epimorphisms. A diagram as in (\ref{regular po}) is always a pushout; it is called a \emph{regular pushout} \cite{B} (alternatively, a \emph{double extension} as in \cite{GRossi}, that was inspired by \cite{Jan})
when, moreover, the canonical morphism $\langle g,c \rangle\colon C\twoheadrightarrow D\times_B A$ to the pullback $D\times_B A$ of $d$ and $f$ is a regular epimorphism.

\begin{remark}
\label{calculus}
{The condition of being a regular pushout can be expressed in terms of the calculus of relations: a commutative diagram of type (\ref{regular po}) is a regular pushout if and only if $c g^{\circ}= f^{\circ}  d$ or, equivalently, $g c^{\circ} = d^{\circ} f$. Observe also that the vertical morphisms $g$ and $f$ are split epimorphisms, so that they induce a split epimorphism from $R_c$ to $R_d$. Consequently, the image of $R_c$ along $g$ is $R_d$, which can be written as: $g\langle R_c \rangle=R_d$, i.e. $gc^{\circ}cg^{\circ}=d^{\circ}d$.}
\end{remark}

The regular categories that are Mal'tsev categories can be characterised through regular pushouts, as it follows from the results of D. Bourn in \cite{B}. For the reader's convenience we give a new simple proof of this fact that uses the calculus of relations, and is also suitable to be extended to the ``relative context'' (see the last section).

\begin{proposition}
\label{Mal'tsev=regular po}
A regular category $\CC$ is a Mal'tsev category if and only if any commutative diagram of the form \emph{(\ref{regular po})} is a regular pushout; equivalently: $$f^{\circ} d = cg^{\circ}.$$
\end{proposition}
\proof
Suppose that $\CC$ is a regular Mal'tsev category. Then:

$$\begin{tabular}{lll}
    $f^{\circ} d$ & $= cc^{\circ}f^{\circ}d$ &  ($c$ is a regular epimorphism; Lemma~\ref{pps of ms as relations}(2)) \\
    & $= cg^{\circ}d^{\circ}d$ &  ($f\cdot c=d \cdot g$) \\
    & $= cg^{\circ} g c^{\circ}cg^{\circ}$ & ($R_d=g\langle R_c \rangle$; see Remark~\ref{calculus}) \\
    & $= cc^{\circ}cg^{\circ}gg^{\circ}$ & ($R_g R_c=R_c R_g$, by Proposition \ref{Definition: Mal'tsev category}) \\
    & $=cg^{\circ}$. & (Lemma~\ref{pps of ms as relations}(1))
\end{tabular}$$

\noindent Conversely, let us consider regular epimorphisms $f\colon X \twoheadrightarrow Y$ and $g\colon X\twoheadrightarrow Z$. We want to prove that $R_f R_g = R_g R_f$. For this we build the following diagram
$$
\xymatrix{
    R_f \ar@<-3pt>[d]_-{f_1} \ar@<3pt>[d]^-{f_2} \ar@{>>}[r]^-{\gamma} & g\langle R_f \rangle \ar@<-3pt>[d]_-{r_1} \ar@<3pt>[d]^-{r_2} \\
    X \ar@{>>}[d]_-f \ar@{>>}[r]_-g & Z \\ Y }
$$
that represents the regular image $g\langle R_f \rangle$ of $R_f$ along $g$. It is easy to see that $g\langle R_f \rangle$ is a reflexive relation, so that we get two commutative diagrams of type (\ref{regular po}). By assumption one then has the equalities: (A) $f_2 \gamma^{\circ}=g^{\circ}r_2$ and (B) $\gamma f_1^{\circ}=
r_1^{\circ}g$. Accordingly:
$$\begin{tabular}[b]{lll}
    $R_f R_g$ &= $f_2f_1^{\circ}g^{\circ}g$ \\
    & $=f_2\gamma^{\circ}r_1^{\circ}g$ & ($g\cdot f_1=r_1\cdot \gamma$) \\
    & $=g^{\circ}r_2r_1^{\circ}g$ & (A) \\
    & $=g^{\circ}r_2\gamma f_1^{\circ}$ & (B) \\
    & $=g^{\circ}gf_2f_1^{\circ}$ & ($r_2\cdot \gamma=g\cdot f_2$) \\
    & $=R_g R_f$.
\end{tabular}$$

\endproof

As first observed by D. Bourn \cite{B}, regular Mal'tsev categories have a strong stability property for regular epimorphisms. Again, we give a new simple proof here to make the paper self-contained.

\begin{lemma}
\label{2.5.7 BB} 
Let $\CC$ be a regular Mal'tsev category. Consider a commutative diagram
\begin{equation}\label{cube}
\xymatrix@C=30pt@R=20pt{
    W \times_D C \ar@<-2pt>[dd]_-k \ar[dr]^-{\gamma} \ar@{.>}[rr]^-v & & Y \times_B A \ar@<-2pt>@{-->}[dd]_(.3){h} \ar[dr]^-{\alpha} \\
    & C \ar@<-2pt>[dd]_(.7){g} \ar@{>>}[rr]^(.2){c} & & A \ar@<-2pt>[dd]_-{f} \\
    W \ar@<-2pt>[uu]_-j \ar[dr]_-{\delta} \ar@{-->>}[rr]^(.7){w} & & Y \ar@<-2pt>@{-->}[uu]_(.7)i \ar@{-->}[dr]_(.4){\beta} \\
    & D \ar@<-2pt>[uu]_(.3)t \ar@{>>}[rr]_-d & & B, \ar@<-2pt>[uu]_-s }
\end{equation}
where the front square is of the form \emph{(\ref{regular po})}, $\beta \cdot w = d \cdot \delta$,
$w$ is a regular epimorphism, \\$(W \times_D C, k, \gamma)$ and $(Y \times_B A, h, \alpha)$ are pullbacks.
 Then the comparison morphism $v \colon W \times_D C \rightarrow Y \times_B A$ is also a regular epimorphism.
\end{lemma}
\proof
The exterior rectangle in the commutative diagram
$$
\xymatrix{W \times_D C \ar[r]^-{\gamma} \ar@<-2pt>[d]_k & C \ar@<-2pt>[d]_g \ar@{->>}[r]^c  & A\ar@<-2pt>[d]_f \\
W  \ar@<-2pt>[u] \ar[r]_{\delta} & D \ar@<-2pt>[u] \ar@{->>}[r]_d & B \ar@<-2pt>[u]
}
$$
is such that the comparison arrow $\langle k, c \cdot \gamma \rangle \colon W \times_D C \rightarrow W \times_B A$ to the pullback of $d \cdot \delta$ along $f$ is a regular epimorphism, as a ``composite'' of a pullback with a regular pushout.
Accordingly, the exterior rectangle in the commutative diagram
$$
\xymatrix{W \times_D C \ar[r]^-{v} \ar@<-2pt>[d]_k & Y \times_B A \ar@<-2pt>[d]_h \ar[r]^-{\alpha}  & A \ar@<-2pt>[d]_f \\
W  \ar@<-2pt>[u] \ar@{->>}[r]_{w} & Y \ar@<-2pt>[u] \ar[r]_{\beta} & B. \ar@<-2pt>[u]
}
$$
has the property that $\langle k, \alpha \cdot v \rangle \colon W \times_D C \rightarrow  W \times_B A$ is a regular epimorphism.
From the fact that the right hand square is a pullback and that $w$ is a regular epimorphism it then easily follows that the arrow $v$ is a regular epimorphism, as desired.
\endproof

We now show that the property considered in Lemma \ref{2.5.7 BB} actually characterises the regular categories which are Mal'tsev categories.

\begin{proposition}\label{thm cube}
Let $\CC$ be a regular category. The following conditions are equivalent:
\begin{enumerate}
\item[(a)] $\CC$ is a Mal'tsev category;
\item[(b)] for any commutative cube \emph{(\ref{cube})}, the comparison morphism \\ $v \colon W \times _D C \rightarrow Y \times_B A$ is a regular epimorphism.
\end{enumerate}
\end{proposition}
\proof
(a) $\Rightarrow$ (b) by Lemma~\ref{2.5.7 BB}.\\
(b) $\Rightarrow$ (a) Consider a commutative diagram of the type (\ref{regular po}). It induces the commutative cube
$$
\xymatrix@C=30pt@R=20pt{
    C \ar@<-2pt>[dd]_-g \ar@{=}[dr] \ar@{.>}[rr]^-v & & D \times_B A \ar@<-2pt>@{-->}[dd]_(.3){h} \ar[dr]^-{\alpha} \\
    & C \ar@<-2pt>[dd]_(.7){g} \ar@{>>}[rr]^(.2){c} & & A \ar@<-2pt>[dd]_-{f} \\
    D \ar@<-2pt>[uu]_-t \ar@{=}[dr] \ar@{=}[rr] & & D \ar@<-2pt>@{-->}[uu]_(.7)i \ar[dr]_(.4){d} \\
    & D \ar@<-2pt>[uu]_(.3)t \ar@{>>}[rr]_-d & & B, \ar@<-2pt>[uu]_-s }
$$
where the right hand face is the pullback of $d$ and $f$. By assumption, $v = \langle g, c \rangle$ is a regular epimorphism, and then $\CC$ is a Mal'tsev category by Proposition~\ref{Mal'tsev=regular po}.
\endproof

\section{The Cuboid Lemma}
\label{The Cuboid Lemma}
In this section we use the characterisation of Mal'tsev categories given in Proposition~\ref{thm cube} to explore a stronger form of the (denormalised) $3 \times 3$ Lemma \cite{B,L}. This follows the same line of research as in \cite{GR}, where the so-called \emph{Goursat pushouts} were used to show that a regular category is a Goursat category if and only if the $3 \times 3$ Lemma holds. 

Recall that the $3 \times 3$ Lemma states the following: given a commutative diagram in $\CC$ of the form
\begin{equation}
\label{3x3 diagram}
\xymatrix@=30pt{
  R_{\bar{g}} \ar@<-3pt>[d]_-{\bar{g}_1} \ar@<3pt>[d]^-{\bar{g}_2} \ar@<-3pt>[r]_-{t_2} \ar@<3pt>[r]^-{ t_1} &
  R_g \ar@<-3pt>[d]_-{g_1} \ar@<3pt>[d]^-{g_2} \ar[r]^-{v} &
  R_f \ar@<-3pt>[d]_-{f_1} \ar@<3pt>[d]^-{f_2} \\
  R_c \ar@<-3pt>[r]_-{c_2} \ar@<3pt>[r]^-{c_1} \ar@{>>}[d]_-{\bar{g}} &
  A \ar@{>>}[d]_-g \ar@{>>}[r]^-c & C \ar@{>>}[d]^-f \\
  S \ar@<-3pt>[r]_-{s_2} \ar@<3pt>[r]^-{s_1}& B \ar[r]_-d & D }
\end{equation}
where the columns are exact (forks), as well as the middle row, then the upper row is exact if and only if the lower row is exact.

Following the analogous terminology introduced in \cite{ZJan10}, the $3 \times 3$ Lemma can be decomposed into two weaker formulations: the \emph{Upper $3 \times 3$ Lemma}, stating that the exactness of the lower row implies the exactness of the upper row, and the \emph{Lower $3 \times 3$ Lemma},  stating that the exactness of the upper row implies the exactness of the lower row.
 It was shown in \cite{GR} (Proposition 1), that these two (apparently) weaker formulations are equivalent to each other, and they are also equivalent to the $3 \times 3$ Lemma: all these formulations provide characterisations of regular Goursat categories (we refer the interested reader also to \cite{Star3x3} for further developments).

By replacing the kernel pairs in the three exact columns of (\ref{3x3 diagram}) with pullbacks of regular epimorphisms along arbitrary morphisms, we get a ``three dimensional version'' of the $3 \times 3$ Lemma. Since these vertical exact forks are replaced by squares, horizontally we should now take $4$ ``forks'' (instead of $3$), two of which are then required to be exact forks. This gives rise to a diagram whose external part is a \emph{cuboid}, and this explains the terminology adopted for the following property.
\vspace{2mm}

\noindent \textbf{The Upper Cuboid Lemma} \\
Let $\CC$ be a regular category, and consider any commutative diagram in $\mathbb C$
\begin{equation}
\label{regular 4x3}
\xymatrix@C=18pt@R=20pt{
    & T \ar@{>>}[dl]_(.6){\bar{k}} \ar@<3pt>[rrrrr]^(.51){t_1} \ar@<-3pt>[rrrrr]_(.51){t_2} \ar[ddr]^(.75){\bar{\gamma}} & & & & &
    V \ar@{-->>}[dl]_(.6){k} \ar[rrrrr]^-v \ar[ddr]^(.75){\gamma} & & & & & X \ar@{-->>}[dl]_(.6){h} \ar[ddr]^(.75){\alpha} \\
    R_w \ar@<3pt>@{-->}[rrrrr]^(.74){w_1} \ar@<-3pt>@{-->}[rrrrr]_(.74){w_2} \ar[ddr]_-(.2){\bar{\delta}} & & & & &
    W \ar@{-->>}[rrrrr]^(.73){w} \ar@{-->}[ddr]_-(.2){\delta} & & & & &
    Y \ar@{-->}[ddr]_-(.2){\beta} \\
    & & R_c \ar@{>>}[dl]_(.6){\bar{g}} \ar@<3pt>[rrrrr]^(.28){c_1} \ar@<-3pt>[rrrrr]_(.28){c_2} & & & & &
    C \ar@{>>}[dl]_(.6){g} \ar@{>>}[rrrrr]_(.27){c} & & & & & A \ar@{>>}[dl]_(.6){f}  \\
    & S \ar@<3pt>[rrrrr]^(.52){s_1} \ar@<-3pt>[rrrrr]_(.52){s_2} & & & & &
    D \ar[rrrrr]_-{d} & & & & & B, }
\end{equation}
where the three diamonds are pullbacks of regular epimorphisms along arbitrary morphisms and the two middle rows are exact forks. Then the Upper Cuboid Lemma holds true in $\mathbb C$ if, for any commutative diagram (\ref{regular 4x3}), \emph{the upper row is exact whenever the lower row is exact}.

\begin{remark}
{Of course, the validity of the Upper Cuboid Lemma is stronger than that of the Upper $3 \times 3$ Lemma (which is itself equivalent to the $3 \times 3$ Lemma, as shown in \cite{GR}). Indeed, the Upper $3 \times 3$ Lemma is simply obtained as the special case of the Upper Cuboid Lemma where one takes the same regular epimorphism twice to build the pullbacks in diagram (\ref{regular 4x3}). As we shall prove in Theorem~\ref{thm upper cuboid}, the validity of the Upper Cuboid Lemma characterises the regular categories that are Mal'tsev categories.}
\end{remark}
We first consider a ``split version'' of the (Upper) Cuboid Lemma to better explore the connection with diagrams of the form (\ref{cube}).
\vspace{2mm}

\noindent\textbf{The Split Cuboid Lemma}\label{CuboidStatement} \\
Let $\CC$ be a regular category, and consider a commutative diagram in $\mathbb C$
\begin{equation}
\label{4x3}
\xymatrix@C=18pt@R=20pt{
    & T \ar@<-2pt>[dl]_(.6){\bar{k}} \ar@<3pt>[rrrrr]^(.51){t_1} \ar@<-3pt>[rrrrr]_(.51){t_2} \ar[ddr]^(.75){\bar{\gamma}} & & & & &
    V \ar@<-2pt>@{-->}[dl]_(.6){k} \ar[rrrrr]^-v \ar[ddr]^(.75){\gamma} & & & & & X \ar@<-2pt>@{-->}[dl]_(.6){h} \ar[ddr]^(.75){\alpha} \\
    R_w \ar@<-2pt>[ur]_-{\bar{j}} \ar@<3pt>@{-->}[rrrrr]^(.74){w_1} \ar@<-3pt>@{-->}[rrrrr]_(.74){w_2} \ar[ddr]_-(.2){\bar{\delta}} & & & & &
    W \ar@<-2pt>@{-->}[ur]_-{j} \ar@{-->>}[rrrrr]^(.73){w} \ar@{-->}[ddr]_-(.2){\delta} & & & & &
    Y \ar@<-2pt>@{-->}[ur]_-{i} \ar@{-->}[ddr]_-(.2){\beta} \\
    & & R_c \ar@<-2pt>[dl]_(.6){\bar{g}} \ar@<3pt>[rrrrr]^(.28){c_1} \ar@<-3pt>[rrrrr]_(.28){c_2} & & & & &
    C \ar@<-2pt>[dl]_(.6){g} \ar@{>>}[rrrrr]_(.27){c} & & & & & A \ar@<-2pt>[dl]_(.6){f}  \\
    & S \ar@<-2pt>[ur]_(.55){\bar{t}} \ar@<3pt>[rrrrr]^(.52){s_1} \ar@<-3pt>[rrrrr]_(.52){s_2} & & & & &
    D \ar@<-2pt>[ur]_-{t} \ar[rrrrr]_-{d} & & & & & B, \ar@<-2pt>[ur]_-{s}}
\end{equation}
where the three diamonds are pullbacks of split epimorphisms along arbitrary morphisms and the two middle rows are exact forks. The Split Cuboid Lemma holds true in $\mathbb C$ if, for any commutative diagram (\ref{4x3}), the \emph{the upper row is exact if and only if the lower row is exact}.


\label{suf conds}
Note that $t_1, t_2, s_1$ and $s_2$ are just parallel morphisms with no common splitting initially required. However, by the commutativity of the diagram, $s_1$ and $s_2$ do actually have a common splitting $\bar{g}\cdot \langle 1_C,1_C \rangle \cdot t$.
The commutativity of the diagram also implies that $d$ is a regular epimorphism, that $d\cdot s_1=d\cdot s_2$ since $\bar{g}$ is an epimorphism and, moreover, that $v\cdot t_1=v\cdot t_2$, because the pair of morphisms $(h,\alpha)$ is jointly monomorphic.

\begin{remark}\label{upper=full}
{Observe that the existence of the split epimorphism $\bar{g}$ implies that we always have $S=R_d$, i.e. the lower row is necessarily exact. Indeed, let $\phi \colon S \rightarrow R_d$ be the canonical arrow induced by the universal property of $R_d$ thanks to the equality $d\cdot s_1=d\cdot s_2$. This arrow $\phi$ is a monomorphism since $S$ is a relation (the pair of arrows $(c_1 \cdot \overline{t}, c_2 \cdot \overline{t})$ is jointly monic by assumption). But $\phi$ is also a split epimorphism, since $\phi \cdot \overline{g}$ is a split epimorphism.
Accordingly, the arrow $\phi \colon S \rightarrow R_d$ is an isomorphism. Consequently, the Split Cuboid Lemma is actually equivalent to the (apparently weaker) \emph{Upper Split Cuboid Lemma}, asserting that the }\emph{upper row in diagram}\emph{ (\ref{4x3})} \emph{is exact whenever the lower row is exact}.
\end{remark}

We are now ready to state our main result:

\begin{theorem}\label{thm upper split cuboid}
Let $\CC$ be a regular category. The following conditions are equivalent:
\begin{enumerate}
\item[(a)] $\CC$ is a Mal'tsev category;
\item[(b)] the Split Cuboid Lemma holds true in $\CC$;
\item[(c)] the Upper Split Cuboid Lemma holds true in $\CC$.
\end{enumerate}
\end{theorem}
\proof
(a) $\Rightarrow$ (b) From the exactness of the lower row in diagram (\ref{4x3}), one clearly has that $T=R_v$, since $S=R_d$ and ``kernel pairs commute with pullbacks''. Consequently, to prove the exactness of the upper row, it suffices to show that $v$ is a regular epimorphism.
The right cube of diagram (\ref{4x3}) is of the form (\ref{cube}), and $v$ is then a regular epimorphism by Proposition~\ref{thm cube}. This shows that the upper row is exact as well.
Conversely, the fact that the lower row is always exact has been proved in Remark \ref{upper=full}.

\noindent (b) $\Rightarrow$ (c) Trivial.

\noindent(c) $\Rightarrow$ (a) Consider a diagram of the form (\ref{cube}) and take the kernel pairs of $c,d$ and $w$. We get an induced split epimorphism $\bar{g}\colon R_c \rightarrow R_d$, with splitting $\bar{t}$, and an induced morphism $\bar{\delta}\colon R_w \rightarrow R_d$. We obtain a diagram of the form (\ref{4x3}) by defining $(T, \bar{k}, \bar{\gamma})$ as the pullback of $\bar{g}$ and $\bar{\delta}$. Applying the Upper Split Cuboid Lemma to this diagram, we conclude that the upper row is exact and, consequently, $v$ is a regular epimorphism. By Proposition~\ref{thm cube}, $\CC$ is a Mal'tsev category.
\endproof

The announced characterisation of Mal'tsev categories in terms of the validity of the Upper Cuboid Lemma can then be given:
\begin{theorem}\label{thm upper cuboid}
Let $\CC$ be a regular category. The following conditions are equivalent:
\begin{enumerate}
\item[(a)] $\CC$ is a Mal'tsev category;
\item[(b)] the Upper Cuboid Lemma holds true in $\CC$.
\end{enumerate}
\end{theorem}
\proof
(a) $\Rightarrow$ (b) Suppose that the lower row is exact. Then $v$ is a regular epimorphism if and only if the following diagram gives the (regular epimorphism, monomorphism) factorisation of the morphism $\langle w\cdot k, c\cdot \gamma \rangle$
$$
\xymatrix@C=20pt{V \ar[r]^-v \ar@{ >->}[d]_-{\langle k, \gamma \rangle} & X \ar@{ >->}[d]^-{\langle h, \alpha \rangle} \\
    W\times C \ar@{>>}[r]_-{w\times c} & Y\times A.}
$$
This translates into having the equality $c\gamma k^{\circ}w^{\circ}=\alpha h^{\circ}$ or, equivalently, the equality $cg^{\circ}\delta w^{\circ}=f^{\circ}\beta$, since the middle and the right diamonds of (\ref{regular 4x3}) are pullbacks (see Remark~\ref{calculus}). This latter can be proved as follows:
$$\begin{tabular}{lll}
    $f^{\circ} \beta$ & $= cc^{\circ}f^{\circ}\beta ww^{\circ}$ &  ($c$, $w$ are a regular epimorphisms; Lemma~\ref{pps of ms as relations}(2)) \\
    & $= cg^{\circ}d^{\circ}d \delta w^{\circ}$ &  ($f\cdot c=d \cdot g$, $\beta\cdot w=d\cdot \delta$) \\
    & $= cg^{\circ} g c^{\circ}cg^{\circ}\delta w^{\circ}$ & ($S=R_d=g\langle R_c \rangle$) \\
    & $= cc^{\circ}cg^{\circ}gg^{\circ}\delta w^{\circ}$ & ($R_g R_c=R_c R_g$, by Proposition \ref{Definition: Mal'tsev category}) \\
    & $=cg^{\circ}\delta w^{\circ}$. & (Lemma~\ref{pps of ms as relations}(1))
\end{tabular}$$
\noindent (b) $\Rightarrow$ (a) This implication is obvious by Theorem~\ref{thm upper split cuboid} and the fact that the validity of the Upper Cuboid Lemma implies the one of the Upper Split Cuboid Lemma.
\endproof

\section{The relative context}
\label{The relative context}

In this section we briefly analyse the so-called ``relative context'' introduced by T. Janelidze ~\cite{TJ}, and verify that our results easily extend from the ``absolute context'' considered in the previous sections to the relative one.

We begin by recalling the notion of a \emph{relative regular category} as defined by J. Goedecke and T. Janelidze \cite{RGC,RS-AC}. As in the previous sections, we shall always assume that the base category is finitely complete.

A \emph{relative regular category} is a pair $(\AA, \EE)$, where $\AA$ is a finitely complete category and $\EE$ is a class of regular epimorphisms in $\AA$ such that:
\begin{itemize}
\item[(E1)]
 $\EE$ contains all isomorphisms;
 \item[(E2)] the pullback of a morphism in $\EE$ belongs to $\EE$;  \item[(E3)] $\EE$ is closed under composition;  \item[(E4)] if $f$ and $g\cdot f\in \EE$, then $g\in \EE$;  \item[(F)\, ] if a morphism factors as $f=e\cdot m$, with $e\in \EE$ and $m$ a monomorphism, then it also factors (essentially uniquely) as $f=m'\cdot e'$, where $m'$ is a monomorphism and $e'\in \EE$.
\end{itemize}
Given a finitely complete category $\AA$ that admits coequalisers of kernel pairs and $\EE$ the class of all regular epimorphisms in $\AA$, then $(\AA,\EE)$ is a relative regular category if and only if $\AA$ is a regular category: this context is usually referred to as the ``absolute context''.

In the extension from the absolute to the relative context one replaces the regular epimorphisms with morphisms that belong to $\EE$.
As shown in \cite{TJ} (see also \cite{RGC}), relative regular categories provide the appropriate setting to develop a suitable calculus of $\EE$-relations. Recall that an \emph{$\EE$-relation} $R$ from $A$ to $B$ is a subobject $\langle r_1,r_2 \rangle\colon R\rightarrowtail A\times B$ such that $r_1, r_2\in \EE$. We can identify a morphism $f:A\rightarrow B$ in $\EE$ with the $\EE$-relation $\langle 1_A,f \rangle\colon A\rightarrowtail A\times B$. The definition of the composition of $\EE$-relations can be given exactly as in the absolute case, and the  properties gathered in Section~\ref{Preliminaries} all extend to $\EE$-relations for relative regular categories. However, for the relative version of Lemma~\ref{pps of ms as relations} we require $f:A\rightarrow B$ to belong to $\EE$: accordingly, the relative version of ~\ref{pps of ms as relations}(2) should simply state that, for any $f\in \EE$, one has the equality $ff^{\circ}=1_B$ (since the arrows in $\EE$ are regular epimorphisms). 
To extend our results to the relative case, we should consider:
\begin{itemize}
    \item[-] \emph{equivalence $\EE$-relations}: $\EE$-relations which are equivalence relations;
    \item[-] \emph{$\EE$-effective equivalence $\EE$-relations}: $\EE$-relations given by kernel pairs of morphisms in $\EE$;
    \item[-] \emph{$\EE$-exact forks}: exact forks whose coequaliser belongs to $\EE$.
\end{itemize}

The \emph{$\EE$-image} of an $\EE$-relation $R$ along a morphism $f\in \EE$ is given by the monomorphism in the $(\EE$, monomorphism)-factorisation of $(f\times f)\cdot \langle r_1,r_2 \rangle$
$$
\xymatrix{
    R \ar[r]^-{\varphi} \ar@{ >->}[d]_-{\langle r_1, r_2 \rangle} & S \ar@{ >->}[d]^-{\langle s_1, s_2 \rangle} \\
    A\times A \ar[r]_-{f\times f} & B\times B,}
$$
which exists by axiom (F) of the definition of a relative regular category; we write $f\langle R \rangle =S$. As in the absolute case, the relation $f\langle R \rangle$ is reflexive (or symmetric) whenever $R$ is reflexive (or symmetric).

Let us then consider a diagram (of type (\ref{regular po}))
\begin{equation}
\label{regular po2}
\xymatrix{
  C \ar@{>>}[r]^-{c} \ar@<-2pt>[d]_-g & A  \ar@<-2pt>[d]_-f \\
  D \ar@{>>}[r]_-d \ar@<-2pt>[u]_-t & B \ar@<-2pt>[u]_-s }
\end{equation}
in a relative regular category, where $f, c, g$ and $d$ now belong to $\EE$. It is called an \emph{$\EE$-regular pushout} when the induced morphism $\langle g,c \rangle\colon C\twoheadrightarrow D\times_B A$ also belongs to $\EE$.
 In the context of a relative regular category $(\AA, \EE)$ the property that any diagram (\ref{regular po2}) in $\AA$ is an $\EE$-regular pushout will be called the \emph{$\EE$-Mal'tsev axiom} in what follows. This axiom has already been studied by T. Everaert, J. Goedecke, T. Janelidze and T. Van der Linden in \cite{EGV,EGJV}, and is similar to the \emph{$\EE$-Goursat axiom} introduced in \cite{RGC} by J. Goedecke and T. Janelidze.

Remark also that the notion of an $\EE$-Mal'tsev category had essentially been introduced in T. Janelidze's PhD doctoral dissertation \cite{TJ} (see Theorem $2.3.6$ therein, for instance), even though the term $\EE$-Mal'tsev had not been explicitly used there. As in the absolute case (Proposition~\ref{Definition: Mal'tsev category}), these categories can be characterised by the commutativity of the composition of $\EE$-effective equivalence $\EE$-relations.

The statements from Remark~\ref{calculus} extend to the relative context by Theorem 2.10 of \cite{RS-AC} and Lemmas 1.9 and 1.11 of \cite{RGC}.
We also obtain relative versions of Propositions ~\ref{Mal'tsev=regular po} and ~\ref{thm cube} with the same proofs, which have been written on purpose in a way that can be now ``understood'' in the relative context:

\begin{proposition}\label{one}
A relative regular category $(\AA, \EE)$ satisfies the $\EE$-Mal'tsev axiom if and only if $$R_f R_g=R_g R_f$$ for any pair of arrows $f$ and $g$ in $\EE$ with the same domain.
\end{proposition}

\begin{proposition}\label{two}
Let $(\AA, \EE)$ be a relative regular category. The following conditions are equivalent:
\begin{enumerate}
\item[(a)] the $\EE$-Mal'tsev axiom holds true in $\AA$;
\item[(b)] $(\AA, \EE)$ has the property that for any commutative cube \emph{(\ref{cube})}, with $f, g, c, d, w\in \EE$, the induced morphism $v \colon W \times _D C \rightarrow Y \times_B A$ belongs to $\EE$.
\end{enumerate}
\end{proposition}

\noindent\textbf{The Relative Cuboid Lemmas}\\
Let $(\AA, \EE)$ be a relative regular category. Consider a commutative diagram~(\ref{4x3}) (resp. ~(\ref{regular 4x3})) where the three diamonds are pullbacks of split (resp. regular) epimorphisms which belong to $\EE$ along arbitrary morphisms and the two middle rows are $\EE$-exact. Then \emph{the upper row is $\EE$-exact if and only if (resp. whenever) the lower row is $\EE$-exact}.

\vspace{3mm}

\noindent With the same proof as in the absolute case, we obtain:

\begin{theorem}\label{three}
Let $(\AA, \EE)$ be a relative regular category. The following conditions are equivalent:
\begin{enumerate}
\item[(a)] the $\EE$-Mal'tsev axiom holds true in $\AA$;
\item[(b)] the Relative (Upper) Split Cuboid Lemma holds true in $\AA$;
\item[(c)] the Relative Upper Cuboid Lemma holds true in $\AA$.
\end{enumerate}
\end{theorem}


\end{document}